\documentclass[10pt]{article}
\usepackage{amsthm, amsmath, amsfonts,	xcolor,	hyperref,	enumerate,	geometry, multirow,	caption, pxfonts} \usepackage{orcidlink}
\newtheorem{theorem}{Theorem}
\newtheorem{example}[theorem]{Example}
\newtheorem{corollary}[theorem]{Corollary}
\newtheorem{definition}[theorem]{Definition}
\newtheorem{proposition}[theorem]{Proposition}
\newtheorem{lemma}[theorem]{Lemma}
\newtheorem{remark}[theorem]{Remark}

\title{\textbf{Riemannian hypermaps to complex space forms}}

\author{Kiran Meena\orcidlink{0000-0002-6959-5853}
, Bayram \c{S}ahin\orcidlink{0000-0002-9372-1151}
	}

\date{}

\begin{document}
	\maketitle
	\begin{abstract}
		\noindent In this paper, we propose the notion of a \textit{hypermap}, defined as a Riemannian map whose range distribution is tangent to a real hypersurface of a complex space form. This construction leads to a new direction in Riemannian map theory. We establish the fundamental Gauss and Codazzi equations for a hypermap and use them to investigate the extrinsic geometry of its range leaves. In particular, when the ambient complex space form is non-Euclidean, we prove that these leaves cannot be totally umbilical and that their shape operator fails to be parallel. We further incorporate the Hopf condition and define \textit{Hopf hypermap}, thereby obtaining the analog of Hopf hypersurfaces. We also examine the umbilical case and show that an umbilical hypermap cannot exist to a non-Euclidean complex space form. This motivates introducing the \textit{$\eta$-umbilical hypermap}, which broadens the class of examples for studying the interaction between the range geometry and the ambient complex structure. Finally, we analyze the principal curvatures associated with these hypermaps, including results on their multiplicities and constancy.
	\end{abstract}
	
	

    \tableofcontents

   \section{Introduction and main results}
    The geometry of real hypersurfaces in complex space forms has become a substantial area of research. It is closely related to several developments in both mathematics and physics. A notable example is Penrose's twistor correspondence \cite{Penrose}, through which projective lines are associated with real hypersurfaces in Minkowski space. This also provides a geometric framework for relating the CR-structure of real hypersurfaces in complex projective spaces to the geometric structures of physical spacetime. In addition, real hypersurfaces occur in a variety of problems involving spectral geometry and geometric analysis. In particular, \cite[Sect.~3.10]{Cecil_Ryan_Book} describes connections with the spectrum of the Laplacian, Willmore and Lagrangian submanifolds, moment maps, the Tits metric, Hamiltonian systems of hydrodynamic type, Hamiltonian stability, wave equations, etc.
    
    In this direction, Takagi's classification \cite{Takagi} of homogeneous real hypersurfaces in complex projective spaces made a major contribution and subsequently motivated numerous investigations into their geometric properties \cite{Kimura, Kimura_Maeda, Okumura, Choe}. One notable feature of the hypersurfaces in Takagi's classification is that their principal curvatures are constant. Furthermore, Berndt \cite{Berndt}, Montiel and Romero \cite{Montiel, Montiel_Romero}, Ki and Suh \cite{Ki_Suh}, Choe \cite{Choe}, and others also obtained classification results for real hypersurfaces in complex hyperbolic spaces. The intrinsic almost contact structure associated with a real hypersurface also leads to distinguished geometric objects, among which the structure vector field and Hopf hypersurfaces play a particularly important role. This has led to extensive studies of Hopf hypersurfaces satisfying additional curvature conditions, including constant principal curvatures in both complex projective and complex hyperbolic spaces \cite{Chen_2001}. These studies appear in the monograph \cite{Cecil_Ryan_Book} and the survey article \cite{NR_1997}. This remains an active research direction. Motivatingly, the authors Meena et al. recently introduced these analogs to the theory of Riemannian submersions in \cite{Meena_Sahin_RHS, MMV}.
    
    On the other hand, Fischer introduced the notion of a Riemannian map (as a generalization of submanifold/isometric immersion and Riemannian submersion) and showed its applications in quantum models, optical physics, various branches of mathematics, etc. \cite{Fischer_1992}. Thereafter, Riemannian maps have been explored with various geometrically structured manifolds and are still an active area of research; see \cite{Akyol_2023, Demir_2026, Lee_2021, Lee_2022, MGSKRB_2026, Meena_Sahin_Shah, Meena_2025, Polat_2025, Polat_Meena, Prasad_2025, Sahin_2025, Sahin_2010, Sahin_2011, Sahin_book, Kiran_JMAA, SMM_GCFI, SMM_Filomat, Yadav_Meena_2023, Jyoti_2024} and references therein. However, despite significant progress, no study in the literature on Riemannian maps considers real hypersurfaces. The main aim of this paper is to introduce and investigate this case. Since this case is analogous to real hypersurfaces, we term it a {\it hypermap} and define it as follows:
    
\begin{definition}\label{def_RHM}
		A Riemannian map $\pi: (M, g_1)$ $\to$ $(N(4c), g_2, J)$ from a Riemannian manifold to a complex space form is said to be a \textit{hypermap} if the leaves of the range space (denoted by $L_\mathcal{R}$) are tangent to real hypersurfaces of $N$.
	\end{definition}

    \begin{remark}
        Throughout this paper, unless otherwise stated, the term \textit{hypermap} refers to a hypermap from a Riemannian manifold to a complex space form.
    \end{remark}

    Let $\pi$ be a Riemannian map from a Riemannian manifold to a K\"ahler manifold $(N, g_2, J)$ such that the leaves of the range space are tangent to real hypersurfaces of $N$ and $\dim(N) = 2n>2$. Let $\xi$ denote a unit normal vector field on the range space. For any vector field $\pi_\ast X \in \Gamma (\mathcal{R})$, we write 
	\begin{equation}\label{cc1}
		J \pi_\ast X = \phi \pi_\ast X + \eta(\pi_\ast X) \xi,~ U = -J \xi,
	\end{equation} 
	where $\phi: \mathcal{R} \to \mathcal{R}$ is an endomorphism such that $\ker (\phi) =  \operatorname{span} \{U\}$, $\eta$ is a $1$-form, $\eta(\pi_\ast X) \xi \in \Gamma(\mathcal{R}^\perp)$, and $U \in {\Gamma(TN)}$ is the \textit{structure vector field}. Then we find
	\begin{equation}\label{cc2}
		\phi^2 \pi_\ast X = -\pi_\ast X + \eta(\pi_\ast X) U, ~\phi U=0, ~\eta(\phi \pi_\ast X)= 0.
	\end{equation}
	Thus, we have the following result:
	\begin{proposition}\label{prop_for_contact_str}
		The leaves of the range space ($L_{\mathcal{R}}$) admit an almost contact structure $\{\phi, \eta, U\}$.
	\end{proposition}

Furthermore, in Sect. \ref{sec_GC_appl}, we derive several fundamental relations, including the Gauss and Codazzi equations. These relations are then used to establish the following results concerning the non-parallelism of the shape operator $\mathcal{S}$ and the nonexistence of umbilical $L_{\mathcal R}$:
\begin{theorem}\label{non-parallelism_of_S}
		Let $\pi: (M, g_1) \to (N(4c), g_2, J)$ be a hypermap. If $c \neq 0$, then 
        \begin{enumerate}[$(1)$]
            \item $\phi \mathcal{S}$ cannot be identically zero;

            \item $\mathcal{S}$ cannot be parallel.
        \end{enumerate}
	\end{theorem}

\begin{theorem}\label{Nonexistence_of_umbilical_LR}
		Let $\pi: (M, g_1) \to (N(4c), g_2, J)$ be a hypermap. If $c \neq 0$, then $\mathcal{S} = \rho I$ cannot hold, even if $\rho$ is not constant. In particular, $L_{\mathcal{R}}$ cannot be totally umbilical.
	\end{theorem}
Furthermore, we introduce the notions of \textit{Hopfness} and \textit{umbilicity} for a hypermap, as stated below: 

\begin{definition}\label{Riem_Hopf_HM}
		A hypermap $\pi:(M,g_1)\to(N(4c),g_2,J)$ is called a \textit{Hopf hypermap} if the structure vector field $U$ is a principal direction of $\mathcal{S}$. That is, there exists a smooth function $\rho$ such that $\mathcal{S} U = \rho U$.
	\end{definition}
    \begin{definition}\label{umbilical_HRM}
		A hypermap $\pi:(M, g_1) \to (N(4c), g_2, J)$ is said to be an \textit{umbilical hypermap} if there exists a differentiable function $\rho$ on $N$ such that
		\begin{equation*}\label{umbilical_RM_def}
			\mathcal{S}_V \pi_\ast X = \rho \pi_\ast X,
		\end{equation*}
		for all $V \in \Gamma(\mathcal{R}^\perp)$ and $\pi_\ast X \in \Gamma(\mathcal{R})$.
	\end{definition}
    
    Within the framework of Hopf (resp. umbilical) hypermaps, all $L_{\mathcal{R}}$ are tangent to Hopf (resp. umbilical) hypersurfaces of the ambient target space. However, by Theorem \ref{Nonexistence_umbilical_RHM}, no umbilical hypermap exists when the target space is a complex non-Euclidean space.  Therefore, we relax it by defining the notion of an \textit{$\eta$-umbilical hypermap} inspired by the definition of a totally $\eta$-umbilical manifold \cite{Kon}, as follows:
    \begin{definition}\label{def_eta_umbilical_RHM}
		A hypermap $\pi:(M, g_1) \to (N(4c\neq 0), g_2, J)$ is said to be an \textit{$\eta$-umbilical hypermap} if there exist differentiable functions $\rho$ and $\sigma$ on $N$ such that
		\begin{equation*}
			\mathcal{S}_V \pi_{\ast} X= \rho \pi_{\ast} X + \sigma \eta(\pi_{\ast}X) U
		\end{equation*}
		for all $V \in \Gamma(\mathcal{R}^\perp)$ and $\pi_\ast X \in \Gamma(\mathcal{R})$. Furthermore, it is also Hopf because $\mathcal{S}_V U = (\rho + \sigma) U$.
	\end{definition}

     Within the framework of $\eta$-umbilical hypermaps, all $L_{\mathcal{R}}$ are tangent to $\eta$-umbilical hypersurfaces of the ambient target space. Under these notions, we characterize the principal curvatures as follows:
     \begin{theorem}\label{cor1_hopf}
		Let $\pi: (M, g_1) \to (N(4c), g_2, J)$ be a hypermap. If $\mathcal{S} \phi = \phi \mathcal{S}$, then $\pi$ is a Hopf hypermap with the principal curvature $\rho$. Further, if $c \neq 0$, then $\rho$ must be constant.
	\end{theorem}

\begin{theorem}\label{thm_three_cpc_of_LR}
		Let $\pi: (M, g_1) \to (N(4c \neq 0), g_2, J)$ be a Hopf hypermap. If $\phi \mathcal{S}_\xi= \mathcal{S}_\xi \phi$, then $L_\mathcal{R}$ have at most three constant principal curvatures. Equivalently, if $\pi$ is a hypermap such that $L_\mathcal{R}$ are of type A, then $L_\mathcal{R}$ have at most three constant principal curvatures.
	\end{theorem}

    \begin{theorem}\label{cor2_hopf}
		Let $\pi: (M, g_1) \to (N(4c), g_2, J)$ be a Hopf hypermap. Let $\rho$ be the principal curvature corresponding to $U$. Then the following assertions hold:
		\begin{enumerate}[$(1)$]
			\item If $c > 0$, then $\rho$ is constant.
			\item If $c=0$, then either $\rho$ is constant or $\mathcal{S}|_{U^\perp}$ has at most two eigenvalues.
		\end{enumerate}
	\end{theorem}

    \begin{theorem}\label{thm2_cpc_of_LR}
		Let $\pi: (M^m, g_1) \to (N^{2n}(4c\neq0), g_2)$ be a Hopf hypermap. If $\mathcal{S}_\xi \phi + \phi \mathcal{S}_\xi = k \phi$ for some constant $k \neq 0$, then any $L_{\mathcal{R}}$ has three constant principal curvatures $\lambda, \mu, \rho$. In addition,
		\begin{enumerate}[$(1)$]
			\item if $\lambda \neq \mu$, then the multiplicities of $\lambda$ and $\mu$ are equal to $n-1$; 
			
			\item if $\lambda = \mu$, then $\pi$ is $\eta$-umbilical.
		\end{enumerate}
	\end{theorem}

    From the previous theorem, we also have the following straightforward result:
    
    \begin{corollary}\label{cor_cpc_of_LR_EURHM}
		Let $\pi: (M, g_1) \to (N(4c \neq 0), g_2, J)$ be an $\eta$-umbilical hypermap. Then $L_\mathcal{R}$ has two constant principal curvatures.
	\end{corollary}

From the preceding main results, we deduce that such Riemannian maps (or hypermaps) have rich and interesting geometry.
	
    \section{Preliminaries}\label{sec_prelims}
	
	In this section, we gather the existing information required for the subsequent sections. \\
	
	Consider a smooth map $\pi: (M, g_1) \to (N, g_2)$ between Riemannian manifolds $M$ and $N$ such that $0 < {\textrm{\normalfont rank }} \pi < \min\{\dim(M), \dim(N)\}$, and its derivative map $\pi_{\ast p}: T_p M \to T_{\pi(p)} N$ at $p$. We denote the kernel space (known as the \textit{vertical space}) of $\pi_\ast$ at $p \in M$ by $\mathcal{V}_p = (\ker \pi_{\ast p})$ and its orthogonal complement space (known as the \textit{horizontal space}) in $T_p M$ by $\mathcal{H}_p = (\ker \pi_{\ast p})^\perp$. Thus, we have the decomposition
	\begin{equation*}
		T_p M = (\ker \pi_{\ast p}) \oplus (\ker \pi_{\ast p})^\perp = \mathcal{V}_p \oplus \mathcal{H}_p.
	\end{equation*}
	Similarly, we also have the decomposition
	\begin{equation*}
		T_{\pi(p)} N = ({\textrm{\normalfont range }} \pi_{\ast p}) \oplus ({\textrm{\normalfont range }} \pi_{\ast p})^\perp = \mathcal{R}_p \oplus \mathcal{R}_p^\perp,
	\end{equation*}
	where $\mathcal{R}_p = ({\textrm{\normalfont range }} \pi_{\ast p})$ denotes the range of $\pi_\ast$ and its orthogonal complementary space in $T_{\pi(p)} N$ is $\mathcal{R}_p^\perp= ({\textrm{\normalfont range }} \pi_{\ast p})^\perp$.
	
	\begin{definition} \cite{Fischer_1992}
		The map $\pi$ is called a \textit{Riemannian map} at each $p \in M$, if $\pi_{\ast p}: \mathcal{H}_p \to \mathcal{R}_p$ is a linear isometry.
	\end{definition}
    
	The map $\pi_\ast$ can be viewed as a section of the bundle $\operatorname{Hom}(TM,\pi^{-1}TN) \rightarrow M$, where $\pi^{-1}TN$ is the pullback bundle. The bundle $\operatorname{Hom}(TM,\pi^{-1}TN)$ has a connection $\nabla$ induced by the Levi-Civita connection $\nabla^M$ and the pullback connection ${\nabla^\pi}$. With this, the symmetric \textit{second fundamental form of $\pi$} is \cite{Nore_1987}
	\begin{equation}\label{2nd_funda_form}
		(\nabla \pi_\ast) (\theta_1,\theta_2) = {\nabla_{\pi_\ast {\theta_1}}^N} \pi_\ast \theta_2 - \pi_\ast({\nabla}_{{\theta_1}}^M \theta_2),
	\end{equation}
	for all $\theta_1, \theta_2 \in \Gamma(TM)$, where $\nabla^N$ is the Levi-Civita connection on $N$. For convenience, we denote the pullback connection by $\nabla^N$. Furthermore, for horizontal vectors $X, Y$, it follows that \cite{Sahin_2010},
	\begin{equation*}
		(\nabla \pi_\ast) (X,Y)\in \Gamma(\mathcal{R}^\perp).
	\end{equation*}
	For any vector field $X$ on $M$ and any section $V$ of $\mathcal{R}^\perp$, we have $\nabla_X^{\pi\bot} V$, which is the orthogonal projection of $\nabla_X^N V$ on $\mathcal{R}^\perp$, where $\nabla^{\pi \bot}$ is a linear connection on $\mathcal{R}^\perp$ such that $\nabla^{\pi \bot} g_2 = 0$. Then, for $\pi$, a bilinear self-adjoint \textit{shape operator} $\mathcal{S}_V$ is defined as \cite[p. 188]{Sahin_book}:
	\begin{equation}\label{shape_operator}
		\nabla_{\pi_\ast X}^N V = -\mathcal{S}_V \pi_\ast X + \nabla_X^{\pi \bot} V.
	\end{equation}
    
	In terms of the second fundamental form and the mean curvature vector field $H$ of the range space, we have the following result:
	\begin{proposition}\label{umbilicity_cond_RM} \cite{Sahin_2011}
		The Riemannian map $\pi: (M, g_1) \to (N, g_2)$ is umbilical if and only if, for all horizontal vectors $X, Y$, we have
		\begin{equation*}
			(\nabla \pi_\ast) (X, Y) = g_1(X, Y)H.
		\end{equation*}
	\end{proposition}
	
	The \textit{Gauss-Codazzi equation} is the fundamental equation in Riemannian geometry that relates the second fundamental form and the curvatures. This equation has been generalized for Riemannian maps as \cite{Sahin_book}
	\begin{align}\label{Gauss_Eqn_RM}
		R^N(\pi_\ast X, \pi_\ast Y) \pi_\ast Z &= -\mathcal{S}_{(\nabla \pi_\ast)(Y, Z)} \pi_\ast X + \mathcal{S}_{(\nabla \pi_\ast)(X, Z)} \pi_\ast Y + \pi_\ast (R^M (X, Y) Z) \nonumber \\&+ (\nabla_X (\nabla \pi_\ast)) (Y, Z) - (\nabla_Y (\nabla \pi_\ast)) (X, Z),
	\end{align}
	where $X, Y, Z \in \Gamma (\mathcal{H})$. In addition, $R^N$ and $R^M$ denote the Riemannian curvature tensors on $N$ and $M$, respectively.
	
	In the following part, we recall the curvature tensor of the complex space form:
	
	\begin{definition} \cite{Yano}
		Let $(N, g_2, J)$ be an almost Hermitian manifold. It is called a {\it K\"ahler manifold} if $(\nabla_{\vartheta_1}^NJ)\vartheta_2=0$ for $\vartheta_1, \vartheta_2 \in \Gamma(TN)$. Furthermore, a K\"ahler manifold having constant holomorphic sectional curvature $4c$ is referred to as a {\it complex space form} and is denoted by $N(4c)$.
	\end{definition}
	
	For the vector fields $\vartheta_1, \vartheta_2, \vartheta_3 \in \Gamma(TN)$, the curvature tensor of $N(4c)$ is given by \cite[p. 80]{Sahin_book}:
	\begin{eqnarray}\label{CSF_Curvature}
		R^N(\vartheta_1,\,\vartheta_2)\vartheta_3 & = & c\,\{ g_2\,(\vartheta_2,\,\vartheta_3)\vartheta_1 - g_2\,(\vartheta_1,\,\vartheta_3)\vartheta_2 + g_2\,(J\vartheta_2,\, \vartheta_3)J\vartheta_1 \nonumber \\&& - g_2\,(J\vartheta_1,\,\vartheta_3)J\vartheta_2 + 2\,g_2\,(\vartheta_1,\,J\vartheta_2)J\vartheta_3 \}.
	\end{eqnarray}

	\section{Gauss-Codazzi equations and applications}\label{sec_GC_appl}

	This section develops several fundamental identities, including the Gauss-Codazzi equations, and examines their geometric consequences. \\
    
	For a hypermap, we have
	\begin{equation}\label{eqn_from_mss}
		\nabla_{\pi_\ast X}^N U = \phi \mathcal{S}\pi_\ast X,
	\end{equation}
	and 
	\begin{equation}\label{B_U_alpha}
		g_2(\nabla_{\pi_\ast X}^N \nabla \rho, \pi_\ast Y) = g_2(\nabla_{\pi_\ast Y}^N \nabla \rho, \pi_\ast X),
	\end{equation}
	where $\rho$ is an eigenvalue of $\mathcal{S}_\xi$.
	
	We also have some additional formulas for hypermaps, as follows:
	
	\begin{proposition}\label{prop_1}
		Let $\pi: (M, g_1) \to (N(4c), g_2, J)$ be a hypermap. Then for $X, Y, Z \in \Gamma(\mathcal{H})$, we have
		\begin{enumerate}[$(1)$]
			\item $\pi_\ast (\nabla_X^M {}^\ast \pi_\ast U) = \phi \mathcal{S}_\xi \pi_\ast X$;
			
			\item $(\nabla_{\pi_\ast X}^N \phi) \pi_\ast Y = \eta (\pi_\ast Y) \mathcal{S} \pi_\ast X - g_2(\mathcal{S} \pi_\ast X, \pi_\ast Y) U + g_2(\pi_\ast Y, \phi \mathcal{S} \pi_\ast X) \xi$;

            \item $g_2(\mathcal{S}\pi_\ast Y, \pi_\ast (\nabla_X^M Z)) = -g_2(\nabla_{\pi_\ast X}^N \mathcal{S} \pi_\ast Y, \pi_\ast Z).$
		\end{enumerate}
	\end{proposition}
	
	\begin{proof}
		Using (\ref{cc1}), we have
		\begin{align*}
			\nabla_{\pi_\ast X}^N U& = - \nabla_{\pi_\ast X}^N J \xi = -J \nabla_{\pi_\ast X}^N \xi.
		\end{align*}
		Since $\dim(\mathcal{R}^\perp)=1$, $\nabla_X^{\pi \perp} \xi = 0$. Then, employing (\ref{shape_operator}), we get $\nabla_{\pi_\ast X}^N U= J \mathcal{S}_\xi \pi_\ast X$. This yields,
		\begin{equation}\label{propeq1}
			\nabla_{\pi_\ast X}^N U = \phi \mathcal{S}_\xi \pi_\ast X + \eta(\mathcal{S}_\xi \pi_\ast X)\xi.
		\end{equation}
		Also,
		\begin{align}\label{propeq2}
			\nabla_{\pi_\ast X}^N U& =\pi_\ast (\nabla_X^M {}^\ast \pi_\ast U) + (\nabla \pi_\ast)(X, {}^\ast \pi_\ast U) \nonumber \\&= \pi_\ast (\nabla_X^M {}^\ast \pi_\ast U)+ \eta(\mathcal{S}_\xi \pi_\ast X)\xi.
		\end{align}
		Then $(1)$ follows from (\ref{propeq1}) and (\ref{propeq2}).
        
  		Further, we have
		\begin{align*}
			(\nabla_{\pi_\ast X}^N \phi) \pi_\ast Y = \nabla_{\pi_\ast X}^N \phi \pi_\ast Y - \phi(\nabla_{\pi_\ast X}^N \pi_\ast Y).
		\end{align*}
		Then using (\ref{cc1}) and (\ref{shape_operator}), we get
		\begin{align}\label{propeq3}
			(\nabla_{\pi_\ast X}^N \phi) \pi_\ast Y & = \nabla_{\pi_\ast X}^N J \pi_\ast Y -\{g_2(\nabla_{\pi_\ast X}^N \pi_\ast Y, U) + g_2(\pi_\ast Y, \nabla_{\pi_\ast X}^N U)\}\xi \nonumber \\&+ \eta(\pi_\ast Y) \mathcal{S}_{\xi} {\pi_\ast X} - \phi(\nabla_{\pi_\ast X}^N \pi_\ast Y).
		\end{align}
		Also, employing (\ref{2nd_funda_form}) and (\ref{cc1}) we have
		\begin{align}\label{propeq4}
			\nabla_{\pi_\ast X}^N J \pi_\ast Y& = J\nabla_{\pi_\ast X}^N \pi_\ast Y = J\pi_\ast (\nabla_X^M Y) + g_2(\mathcal{S}\pi_\ast X, \pi_\ast Y) J\xi \nonumber \\& = \phi (\nabla_{\pi_\ast X}^N \pi_\ast Y) + \eta(\nabla_{\pi_\ast X}^N \pi_\ast Y)\xi - g_2(\mathcal{S}\pi_\ast X, \pi_\ast Y) U.
		\end{align} 
		Using (\ref{propeq4}) in (\ref{propeq3}), we get
		\begin{equation*}
			(\nabla_{\pi_\ast X}^N \phi) \pi_\ast Y = -g_2(\mathcal{S}\pi_\ast X, \pi_\ast Y) U + g_2(\pi_\ast Y, \mathcal{S}_U \pi_\ast X)\xi + \eta(\pi_\ast Y) \mathcal{S}_{\xi} {\pi_\ast X}.
		\end{equation*}
		This implies the proof of $(2)$.
        
        Finally, using (\ref{2nd_funda_form}), we have $g_2(\mathcal{S}\pi_\ast Y, \pi_\ast (\nabla_X^M Z))=g_2(\mathcal{S} \pi_\ast Y, \nabla_{\pi_\ast X}^N \pi_\ast Z)$. Then, using the metric compatibility, we get 
		\begin{equation*}
			g_2(\mathcal{S}\pi_\ast Y, \pi_\ast (\nabla_X^M Z))=\nabla_{\pi_\ast X}^N g_2((\nabla \pi_\ast)(Y, Z), \xi) - g_2(\nabla_{\pi_\ast X}^N \mathcal{S} \pi_\ast Y, \pi_\ast Z).
		\end{equation*}
		Since $\dim (\mathcal{R}^\perp)=1$, we can put $(\nabla \pi_\ast)(Y, Z)= f \xi$, where $f$ is a smooth function on $N$. Hence
		\begin{align*}
			g_2(\mathcal{S}\pi_\ast Y, \pi_\ast (\nabla_X^M Z))=\nabla_{\pi_\ast X}^N (f g_2(\xi, \xi)) - g_2(\nabla_{\pi_\ast X}^N \mathcal{S} \pi_\ast Y, \pi_\ast Z).
		\end{align*}
		Then the required proof of $(3)$ follows because $f$ is real, so its derivative is zero at a point.
	\end{proof}
	Now, we present the Gauss and Codazzi equations in the next two lemmas:
	\begin{lemma}\label{Gauss_Eqn_RHM}
		Let $\pi: (M, g_1) \to (N(4c), g_2, J)$ be a hypermap. Then for $X, Y, Z \in \Gamma(\mathcal{H})$, we have the Gauss equation
		\begin{equation*}
			\pi_\ast (R^M(X, Y)Z) = c\{(\pi_\ast X \wedge \pi_\ast Y)\pi_\ast Z + (\phi \pi_\ast X \wedge \phi \pi_\ast Y)\pi_\ast Z + 2 g_2(\pi_\ast X, \phi \pi_\ast Y) \phi \pi_\ast Z\}.
		\end{equation*}
	\end{lemma}
	
	\begin{proof}
		Since $N$ is a complex space form, (\ref{Gauss_Eqn_RM}) and (\ref{CSF_Curvature}) yield
		\begin{align*}
			g_2(\pi_\ast(R^M(X,Y)Z), \pi_\ast T) =& c\{g_2(\pi_\ast Y, \pi_\ast Z) g_2(\pi_\ast X, \pi_\ast T) - g_2(\pi_\ast X, \pi_\ast Z) g_2(\pi_\ast Y, \pi_\ast T)\\&+g_2(J \pi_\ast Y, \pi_\ast Z) g_2(J \pi_\ast X, \pi_\ast T) - g_2(J \pi_\ast X, \pi_\ast Z) g_2(J \pi_\ast Y, \pi_\ast T)\\&+2g_2(\pi_\ast X, J \pi_\ast Y) g_2(J \pi_\ast Z, \pi_\ast T)\},
		\end{align*}
		where $X,Y,Z,T \in \Gamma(\mathcal{H})$. Using (\ref{cc1}) in the aforementioned equation, we obtain
		\begin{align*}
			\pi_\ast(R^M(X,Y)Z) =& c\{g_2(\pi_\ast Y, \pi_\ast Z) \pi_\ast X- g_2(\pi_\ast X, \pi_\ast Z) \pi_\ast Y\\&+g_2(\phi \pi_\ast Y, \pi_\ast Z) \phi \pi_\ast X - g_2(\phi \pi_\ast X, \pi_\ast Z) \phi \pi_\ast Y\\&+2g_2(\pi_\ast X, \phi \pi_\ast Y) \phi \pi_\ast Z\}.
		\end{align*}
		This implies the required proof.
	\end{proof}
    
	\begin{lemma}\label{Codazzi_Eqn_RHM}
		Let $\pi: (M, g_1) \to (N(4c), g_2, J)$ be a hypermap. Then for $X, Y \in \Gamma(\mathcal{H})$, we have the Codazzi equation
		\begin{align*}
			(\nabla_{\pi_\ast X}^N \mathcal{S}) \pi_{\ast} Y - (\nabla_{\pi_\ast Y}^N \mathcal{S}) \pi_{\ast} X = c \{g_2(\pi_\ast X, U) \phi \pi_{\ast} Y - g_2(\pi_\ast Y, U) \phi \pi_{\ast} X + 2 g_2(\pi_{\ast} X, \phi \pi_\ast Y)U\},
		\end{align*}
		where $U$ is the structure vector field.
	\end{lemma}
	
	\begin{proof}
		Since $\dim(\mathcal{R}^\perp) = 1$, $\nabla_{X}^{\pi \perp} \xi = 0$ and $\nabla_{X}^{\pi \perp} ((\nabla \pi_\ast)(X, Z)) = 0$. With these facts,
		\begin{align*}
			(\nabla_X(\nabla \pi_\ast)) (Y, Z)& = -(\nabla \pi_\ast) (\nabla_X^M Y, Z) - (\nabla \pi_\ast) (Y, \nabla_X^M Z)\\&=-g_2(\mathcal{S} \pi_\ast (\nabla_X^M Z), \pi_\ast Z) \xi - g_2(\mathcal{S} \pi_\ast Y, \pi_\ast (\nabla_X^M Z))\xi.
		\end{align*}
		Employing Proposition \ref{prop_1}, the aforementioned equation yields
		\begin{equation*}
			(\nabla_X(\nabla \pi_\ast)) (Y, Z)= g_2(\nabla_{\pi_\ast X}^N \mathcal{S} \pi_\ast Y - \mathcal{S} (\nabla_{\pi_\ast X}^N \pi_\ast Y), \pi_\ast Z) \xi.
		\end{equation*}
		Then the definition of the covariant derivative yields
		\begin{equation}\label{eq10}
			(\nabla_X(\nabla \pi_\ast)) (Y, Z)= g_2((\nabla_{\pi_\ast X}^N \mathcal{S})\pi_\ast Y, \pi_\ast Z)\xi.
		\end{equation}
		Similarly, we get
		\begin{equation}\label{eq11}
			(\nabla_Y (\nabla \pi_\ast)) (X, Z)= g_2((\nabla_{\pi_\ast Y}^N \mathcal{S})\pi_\ast X, \pi_\ast Z)\xi.
		\end{equation}
		Thus, by (\ref{eq10}) and (\ref{eq11}), we have
		\begin{align}\label{eq12}
			g_2((\nabla_X(\nabla \pi_\ast)) (Y, Z), \xi) - g_2((\nabla_Y (\nabla \pi_\ast)) (X, Z), \xi) = g_2((\nabla_{\pi_\ast X}^N \mathcal{S})\pi_\ast Y - (\nabla_{\pi_\ast Y}^N \mathcal{S})\pi_\ast X, \pi_\ast Z).
		\end{align}
		Since $N$ is a complex space form, by (\ref{Gauss_Eqn_RM}) we have,
		\begin{align*}
			&g_2((\nabla_X \nabla \pi_\ast)(Y, Z), \xi) - g_2((\nabla_Y \nabla \pi_\ast)(X, Z), \xi) = c\{g_2(J \pi_\ast Y, \pi_\ast Z) g_2(\pi_\ast X, -J \xi) \\&-g_2(J \pi_\ast X, \pi_\ast Z) g_2(\pi_\ast Y, -J \xi) + 2g_2(\pi_\ast X, J \pi_\ast Y) g_2(\pi_\ast Z, -J \xi)\}.
		\end{align*}
		Then, using (\ref{eq12}) and (\ref{cc1}) in the aforementioned equation, we obtain
		\begin{align}\label{Codazzi_RHM}
			&g_2((\nabla_{\pi_\ast X}^N \mathcal{S})\pi_\ast Y - (\nabla_{\pi_\ast Y}^N \mathcal{S})\pi_\ast X, \pi_\ast Z)= c\{g_2(\phi \pi_\ast Y, \pi_\ast Z) g_2(\pi_\ast X, U) \nonumber \\&-g_2(\phi \pi_\ast X, \pi_\ast Z) g_2(\pi_\ast Y, U) + 2g_2(\pi_\ast X, \phi \pi_\ast Y) g_2(\pi_\ast Z, U)\}.
		\end{align}
		Then the required proof follows by comparing the two sides.
	\end{proof}
	
	\begin{corollary}\label{corollary1}
		Under the hypothesis of the hypermap, we have
		\begin{enumerate}[$(1)$]
			\item $g_2((\nabla_{\pi_\ast X}^N \mathcal{S}) \pi_\ast Y - (\nabla_{\pi_\ast Y}^N \mathcal{S})\pi_\ast X, U) = 2c g_2 (\pi_\ast X, \phi \pi_\ast Y)$;
			
			\item $g_2((\nabla_{\pi_\ast X}^N \mathcal{S})U, U) = g_2((\nabla_{U}^N \mathcal{S})\pi_\ast X, U) = g_2((\nabla_{U}^N \mathcal{S})U, \pi_\ast X)$.
		\end{enumerate}
	\end{corollary}
    
	\begin{proof}
		The proof follows easily from the previous lemma. 
	\end{proof}
	
	We now prove Theorem $\ref{non-parallelism_of_S}$.
	
	\begin{proof}[Proof of Theorem $\ref{non-parallelism_of_S}$]
		By Proposition \ref{prop_1}, $\pi_\ast (\nabla_X^M {}^\ast \pi_\ast U) = 0$ if and only if $\phi \mathcal{S}_\xi \pi_\ast X = 0$. That means $\mathcal{S}_\xi \pi_\ast X \in \operatorname{span}\{U\}$. Suppose that this condition holds for all $\pi_\ast X \in \Gamma(\mathcal{R})$, then
		\begin{equation}\label{cor_eqn_3}
			\mathcal{S} \pi_\ast X = g_2(\mathcal{S} \pi_\ast X, U)U.
		\end{equation}
		Also, we have
		\begin{equation*}
			(\nabla_{\pi_\ast X}^N \mathcal{S}) \pi_\ast Y = \nabla_{\pi_\ast X}^N \mathcal{S} \pi_\ast Y - \mathcal{S}(\nabla_{\pi_\ast X}^N \pi_\ast Y)
		\end{equation*}
		for all $\pi_\ast X, \pi_\ast Y \in \Gamma(\mathcal{R})$. Applying $\phi$, and then using (\ref{cor_eqn_3}), we get
		\begin{equation*}
			\phi((\nabla_{\pi_\ast X}^N \mathcal{S}) \pi_\ast Y) = 0.
		\end{equation*}
		Hence $(\nabla_{\pi_\ast X}^N \mathcal{S}) \pi_\ast Y \in \operatorname{span}\{U\}$. Then, by the Codazzi equation, we obtain
		\begin{equation*}
			\phi c \{g_2(\pi_\ast X, U) \phi \pi_\ast Y - g_2(\pi_\ast Y, U) \phi \pi_\ast X\} = 0.
		\end{equation*}
		Taking $\pi_\ast Y = U$, we get $-c \phi \pi_\ast X \in \operatorname{span}\{U\}$. Then by (\ref{cor_eqn_3}), we write
		\begin{equation*}
			-c \phi \pi_\ast X = g_2(c \phi \pi_\ast X, U)U = c g_2(\phi \pi_\ast X, U)U = c g_2(\pi_\ast X, \phi U)U = 0.
		\end{equation*}
		Thus, $c=0$ or $\phi \pi_\ast X =0$. If $c \neq 0$, then $\phi \pi_\ast X = 0$ for all $\pi_\ast X \in \Gamma(\mathcal{R})$. This implies $\text{rank}~ \phi \neq 2n-2$. Thus, we get a contradiction. Hence, the first statement follows. Further, suppose $\mathcal{S}$ is parallel, that is, $(\nabla^N \mathcal{S}) = 0$. Then the Codazzi equation gives
		\begin{equation*}
			c \{g_2(\pi_\ast X, U) \phi \pi_\ast Y - g_2(\pi_\ast Y, U)\phi \pi_\ast X + 2 g_2(\pi_\ast X, \phi \pi_\ast Y) U\} = 0.
		\end{equation*}
        Applying $\phi$ on both sides, we get 
		\begin{equation*}
			c \{g_2(\pi_\ast X, U) \phi \pi_\ast Y - g_2(\pi_\ast Y, U)\phi \pi_\ast X\} \in \operatorname{span}\{U\}.
		\end{equation*}
		Taking $\pi_\ast Y = U$, we get $c \phi \pi_\ast X = 0$. Since $c \neq 0$, $\phi \pi_\ast X = 0$ for any $\pi_\ast X \in \Gamma(\mathcal{R})$. This implies a contradiction. Hence $\mathcal{S}$ cannot be parallel.
	\end{proof}
	
	Finally, motivated by \cite{Tashiro_Tachibana}, we prove the nonexistence of the totally umbilical $L_\mathcal{R}$ (Theorem $\ref{Nonexistence_of_umbilical_LR}$) as follows:
	
	\begin{proof}[Proof of Theorem $\ref{Nonexistence_of_umbilical_LR}$]
		The definition of covariant derivative and $\mathcal{S} = \rho I$ yield
		\begin{equation*}
			(\nabla_{\pi_\ast X}^N \mathcal{S}) \pi_\ast Y = \pi_\ast X (\rho)\pi_\ast Y,
		\end{equation*}
		and
		\begin{equation*}
			(\nabla_{\pi_\ast Y}^N \mathcal{S}) \pi_\ast X = \pi_\ast Y (\rho)\pi_\ast X.
		\end{equation*}
		Using the above equations in the Codazzi equation, we obtain
		\begin{align*}
			(\pi_\ast X (\rho)) \pi_\ast Y - (\pi_\ast Y (\rho)) \pi_\ast X = c \{g_2(\pi_\ast X, U) \phi \pi_\ast Y - g_2(\pi_\ast Y, U)\phi \pi_\ast X + 2 g_2(\pi_\ast X, \phi \pi_\ast Y) U\}.
		\end{align*}
		Putting $\pi_\ast Y = U$ in the aforementioned equation, we get
		\begin{equation*}
			\pi_\ast X (\rho) U - U (\rho) \pi_\ast X = -c \phi \pi_\ast X.
		\end{equation*}
        Taking the inner product with $\phi \pi_\ast X$ gives $c=0$. Thus, we get a contradiction. Hence, the proof follows.
	\end{proof}

	\section{Some intermediate results}\label{sec_Hopf_case}
    
	This section presents the intermediate results needed to study the behavior of the principal curvatures. We prove them one by one as follows:
	\begin{lemma}\label{Hopf_thm_1}
		Let $\pi: (M, g_1) \to (N(4c), g_2, J)$ be a Hopf hypermap. Then we have
		\begin{equation*}
			\mathcal{S} \phi \mathcal{S} - \frac{\rho}{2}(\mathcal{S} \phi + \phi \mathcal{S}) - c \phi = 0.
		\end{equation*}
	\end{lemma}
	\begin{proof}
		By using (\ref{eqn_from_mss}), we have
		\begin{equation*}
			(\nabla_{\pi_\ast X}^N \mathcal{S}) U = \nabla_{\pi_\ast X}^N \mathcal{S} U - \mathcal{S}(\phi \mathcal{S} \pi_\ast X).
		\end{equation*}
		Since $\pi$ is a Hopf hypermap, $\mathcal{S}U = \rho U$. Then, using (\ref{eqn_from_mss}) again, the aforementioned equation yields 
		\begin{equation}\label{eq3.13}
			(\nabla_{\pi_\ast X}^N \mathcal{S}) U= \pi_\ast X (\rho) U + (\rho I -\mathcal{S}) \phi \mathcal{S} \pi_\ast X.
		\end{equation}
		Taking the inner product with $U$, we get
		\begin{align}\label{eq3.14}
			g_2((\nabla_{\pi_\ast X}^N \mathcal{S}) U, U) = \pi_\ast X (\rho).
		\end{align}
Therefore, by Corollary \ref{corollary1}, we have 
		\begin{align}\label{eq3.15}
			g_2((\nabla_{\pi_\ast X}^N \mathcal{S})U, U) &= g_2((\nabla_{U}^N \mathcal{S})U, \pi_\ast X) = \pi_\ast X (\rho) = g_2(\pi_\ast X, \nabla^N \rho) \nonumber \\& = g_2(\pi_\ast X, U(\rho)U) = U(\rho) g_2(\pi_\ast X, U).
		\end{align}
		The Eqs. (\ref{eq3.13}), (\ref{eq3.14}), and (\ref{eq3.15}) yield
		\begin{equation*}
			(\nabla_{\pi_\ast X}^N \mathcal{S}) U = U(\rho) g_2(\pi_\ast X, U)U + (\rho I -\mathcal{S}) \phi \mathcal{S} \pi_\ast X.
		\end{equation*}
Taking the inner product with $\pi_\ast Y$ and then using Corollary \ref{corollary1}, we obtain
		\begin{align}\label{eq3.16}
			&g_2((\nabla_{\pi_\ast X}^N \mathcal{S}) U, \pi_\ast Y) = g_2((\nabla_{\pi_\ast X}^N \mathcal{S}) \pi_\ast Y, U) \nonumber \\&= U(\rho) g_2(\pi_\ast X, U) g_2(U, \pi_\ast Y) + g_2((\rho I -\mathcal{S}) \phi \mathcal{S} \pi_\ast X, \pi_\ast Y).
		\end{align}
		Interchanging $\pi_\ast X$ and $\pi_\ast Y$, we have
		\begin{equation}\label{eq3.17}
			g_2((\nabla_{\pi_\ast Y}^N \mathcal{S}) \pi_\ast X, U) = U(\rho) g_2(\pi_\ast X, U) g_2(U, \pi_\ast Y) + g_2((\rho I -\mathcal{S}) \phi \mathcal{S} \pi_\ast Y, \pi_\ast X).
		\end{equation}
Then the difference of (\ref{eq3.16}) and (\ref{eq3.17}) gives
		\begin{equation}\label{eq3.18}
			g_2((\nabla_{\pi_\ast X}^N \mathcal{S}) \pi_\ast Y - (\nabla_{\pi_\ast Y}^N \mathcal{S}) \pi_\ast X, U) = \rho g_2((\mathcal{S} \phi + \phi \mathcal{S}) \pi_\ast X, \pi_\ast Y) - 2g_2(\mathcal{S}\phi \mathcal{S} \pi_\ast X, \pi_\ast Y).
		\end{equation}
		Also, by the Codazzi equation, we have
		\begin{equation}\label{eq3.19}
			g_2((\nabla_{\pi_\ast X}^N \mathcal{S}) \pi_\ast Y - (\nabla_{\pi_\ast Y}^N \mathcal{S}) \pi_\ast X, U) = 2cg_2(\pi_\ast X, \phi \pi_\ast Y).
		\end{equation}
		Then the Eqs. (\ref{eq3.18}) and (\ref{eq3.19}) imply 
		\begin{equation}\label{3.16}
			\rho g_2((\mathcal{S} \phi + \phi \mathcal{S}) \pi_\ast X, \pi_\ast Y) - 2g_2(\mathcal{S}\phi \mathcal{S} \pi_\ast X, \pi_\ast Y) = -2cg_2(\phi \pi_\ast X, \pi_\ast Y).
		\end{equation}
		This implies the required proof.
	\end{proof}
	
	For $U \in \Gamma(\mathcal{R})$, define $U^\perp := \{\pi_\ast X \in \Gamma(TN): g_2(\pi_{\ast} X, U) = 0\}$. Then the following holds.
	
	\begin{corollary}\label{Cor_Hopf}
		Let $\pi: (M, g_1) \to (N(4c), g_2, J)$ be a Hopf hypermap. Then the following statements are valid:
		\begin{enumerate}[$(1)$]
			\item For any $\pi_{\ast} X \in U^\perp$ satisfying $\mathcal{S} \pi_{\ast} X = \lambda \pi_\ast X$, we have
			\begin{equation*}
				\left(2\lambda - \rho\right) \mathcal{S} \phi \pi_{\ast} X - \left(\rho \lambda+2c\right)\phi \pi_{\ast} X = 0.
			\end{equation*}
			
			\item For any $\pi_{\ast} X \in U^\perp$ satisfying $\mathcal{S} \pi_{\ast} X = \lambda \pi_\ast X$ and $\mathcal{S} \phi \pi_\ast X = \mu \phi \pi_\ast X$, we have
			\begin{equation*}
				2 \lambda \mu -{\rho (\lambda + \mu)} -2 c=0.
			\end{equation*}
			
			\item If $T_\lambda := \{W \in \Gamma(TN): \mathcal{S} W = \lambda W\}$ is $\phi$-invariant, i.e. $\phi(T_\lambda) \subset T_\lambda$, then for any $\pi_{\ast} X \in T_\lambda$ satisfying $\mathcal{S} \pi_{\ast} X = \lambda \pi_\ast X$, we have
			\begin{equation*}
				\lambda^2 - \rho \lambda - c = 0.
			\end{equation*}
		\end{enumerate}
	\end{corollary}
	
	\begin{proof}
		Since the proof follows easily from the previous lemma, we omit it.
	\end{proof}
	
	\begin{lemma}\label{thm_hopf_2}
		Let $\pi: (M, g_1) \to (N(4c), g_2, J)$ be a Hopf hypermap. Then we have
		\begin{equation*}
			U(\rho)(\mathcal{S} \phi + \phi \mathcal{S}) = 0.
		\end{equation*}
	\end{lemma}
	
	\begin{proof}
		Let $\sigma = U(\rho)$, so that $\nabla \rho = \sigma U$. Then using (\ref{B_U_alpha}) we get,
		\begin{equation*}
			g_2(\pi_\ast X (\sigma) U + \sigma \nabla_{\pi_\ast X}^N U, \pi_\ast Y) = g_2(\pi_\ast Y (\sigma) U + \sigma \nabla_{\pi_\ast Y}^N U, \pi_\ast X).
		\end{equation*}
		Using (\ref{eqn_from_mss}) in the above equation, we get
		\begin{equation}\label{eq25}
			\pi_\ast X (\sigma) g_2(U, \pi_\ast Y) - \pi_\ast Y (\sigma) g_2(U, \pi_\ast X) + \sigma g_2((\mathcal{S} \phi + \phi \mathcal{S})\pi_\ast X, \pi_\ast Y) = 0.
		\end{equation}
		Setting $\pi_\ast Y = U$ in (\ref{eq25}), we get
		\begin{equation*}
			\pi_\ast X (\sigma) - U (\sigma) g_2(U, \pi_\ast X) + \sigma g_2(\phi \pi_\ast X, \mathcal{S}U).
		\end{equation*}
		Since $\pi$ is a Hopf hypermap, $\mathcal{S} U = \rho U$, where $\rho = g_2(\mathcal{S}U, U)$. Then the above equation yields
		\begin{equation*}
			\pi_\ast X (\sigma) = U(\sigma) g_2(U, \pi_\ast X).
		\end{equation*}
		Similarly, we get
		\begin{equation*}
			\pi_\ast Y (\sigma) = U(\sigma) g_2(U, \pi_\ast Y).
		\end{equation*}
		Then (\ref{eq25}) yields
		\begin{equation*}
			g_2(\sigma(\mathcal{S} \phi + \phi \mathcal{S})\pi_\ast X, \pi_\ast Y) = 0.
		\end{equation*}
		Hence, the proof follows.
	\end{proof}

    \begin{proposition}
		Let $\pi: (M^m, g_1) \to (N^n(4c), g_2, J)$ be a Hopf hypermap. If $\phi \mathcal{S}_\xi + \mathcal{S}_\xi \phi= 0$, then $c \leq 0$.
	\end{proposition}
	
	\begin{proof}
		Since $\phi \mathcal{S}_\xi + \mathcal{S}_\xi \phi = 0$, $\phi \mathcal{S}_\xi U = 0$. This implies $\mathcal{S}_\xi U= \eta(\mathcal{S}_\xi U) U$. That is, $U$ is an eigenvector of $\mathcal{S}_\xi$. Then (\ref{3.16}) implies
		\begin{align*}
			c g_2(\pi_\ast X, \phi \pi_\ast Y) = -g_2(\mathcal{S}_\xi \phi \mathcal{S}_\xi \pi_\ast X, \pi_\ast Y)=g_2(\mathcal{S}_\xi \pi_\ast X, \phi \mathcal{S}_\xi \pi_\ast Y).
		\end{align*}
		Replacement of $\pi_\ast X$ by $\phi \pi_\ast X$ yields 
		\begin{align*}
			c g_2(\phi \pi_\ast X, \phi \pi_\ast Y) = g_2(\mathcal{S}_\xi \phi \pi_\ast X, \phi \mathcal{S}_\xi \pi_\ast Y) = -g_2(\mathcal{S}_\xi \phi \pi_\ast X, \mathcal{S}_\xi \phi \pi_\ast Y).
		\end{align*}
		Now, by taking $\pi_\ast X = \pi_\ast Y$, we obtain
		\begin{equation*}
			c g_2(\phi \pi_\ast X, \phi \pi_\ast X) = -g_2(\mathcal{S}_\xi \phi \pi_\ast X, \mathcal{S}_\xi \phi \pi_\ast X) \leq 0.
		\end{equation*}
		Since $\textrm{\normalfont rank } \phi = 2n-2$ and $n > 1$, we must have $c\leq 0$.
	\end{proof}

    \section{Behavior of principal curvatures}
	
	In this section, we prove the main results (specifically, the constancy and enumeration of the principal curvatures) as follows:
	
	\begin{proof}[Proof of Theorem $\ref{cor1_hopf}$]
		Since $\phi \mathcal{S} = \mathcal{S} \phi$, $\phi \mathcal{S} U = \mathcal{S} \phi U = 0$. That is, $\mathcal{S} U \in \operatorname{span}\{U\} = \ker \phi$. Hence, $\mathcal{S} U = \rho U$. Equivalently, $\pi$ is a Hopf hypermap.
		Then, by the previous lemma, if $\sigma = U (\rho) \neq 0$ at some point $\pi(p)$, we have $\phi \mathcal{S} + \mathcal{S} \phi = 0$. Therefore, $\phi \mathcal{S} = - \mathcal{S} \phi = -\phi \mathcal{S}$. This implies
		\begin{equation*}
			\phi \mathcal{S} = \mathcal{S} \phi = 0.
		\end{equation*}
        Then, using the Lemma \ref{Hopf_thm_1}, we obtain $c \phi \pi_\ast X = 0$. Thus, if $c \neq 0$, then $\phi \pi_\ast X= 0$ for any $\pi_\ast X \in \operatorname{span}\{U\}$. Thus, we reach a contradiction, so $\sigma = U(\rho) = 0$. This completes the required proof.
	\end{proof}
	
	\begin{proof}[Proof of Theorem $\ref{thm_three_cpc_of_LR}$]
		Since $\mathcal{S}_\xi \phi = \phi\mathcal{S}_\xi$, $\phi\mathcal{S}_\xi U= 0$. This yields
		\begin{equation*}
			J\mathcal{S}_\xi U= \eta(\mathcal{S}_\xi U) \xi~\text{or}~\mathcal{S}_\xi U = \eta(\mathcal{S}_\xi U) U.
		\end{equation*}
		That is, $U$ is an eigenvector of $\mathcal{S}_\xi$. Then by (\ref{eq25}) we obtain,
		\begin{equation*}
			\sigma g_2(\phi\mathcal{S}_\xi \pi_\ast X, \pi_\ast Y)= 0.
		\end{equation*}
		If $\sigma \neq 0$ at some point $\pi(p)$, then $\phi \mathcal{S}_\xi \pi_\ast X= 0$ and hence (\ref{3.16}) implies that
		\begin{equation*}
			c g_2(\pi_\ast X, \phi \pi_\ast Y)=0.
		\end{equation*}
		Thus $c=0$, a contradiction since $c \neq 0$. So, $\sigma$ must be zero, and hence $\rho$ is constant.
        
        Also, from (\ref{3.16}) it follows that 
		\begin{equation*}
			\phi\mathcal{S}_\xi^2 \pi_\ast X- \rho \phi \mathcal{S}_\xi \pi_\ast X- c \phi \pi_\ast X= 0.
		\end{equation*}
		Then, using (\ref{cc2}), we get 
		\begin{equation*}
			\mathcal{S}_\xi^2 \pi_\ast X - \rho \mathcal{S}_\xi \pi_\ast X - c \pi_\ast X + c \eta(\pi_\ast X) U = 0.
		\end{equation*} 
		The $\mathcal{S}_\xi$ is a self-adjoint operator, and let $\{\pi_\ast e_i\}_{i=1}^{2n-2}$ be an orthonormal basis of $T_{\pi(p)} L_{\mathcal{R}}$. Furthermore, we may assume
		\begin{equation*}
			\mathcal{S}_\xi \pi_{\ast} e_i = {\lambda}_i \pi_{\ast} e_i ~\text{for}~ i= 1, 2, \dots, 2n-2 ~\text{and}~\mathcal{S}_\xi e_{2n-1} =\mathcal{S}_\xi U = \rho e_{2n-1}.
		\end{equation*}
		Thus, we get at most two distinct eigenvalues $\lambda_i$, denoted by constants $\lambda$ and $\mu$ (say). Also, we have
		\begin{equation*}
			\lambda + \mu = \textrm{trace}~\mathcal{S}_\xi = \rho,
		\end{equation*}
		and
		\begin{equation*}
			\lambda \mu = \det~\mathcal{S}_\xi = -c.
		\end{equation*}
		This proves the result.
	\end{proof}
    
	\begin{proof}[Proof of Theorem $\ref{cor2_hopf}$]
		Since $\pi$ is a Hopf hypermap, $\mathcal{S} U = \rho U$. Then, by Lemma \ref{thm_hopf_2}, if $\sigma = U (\rho) \neq 0$ at some point $\pi(p)$, we have $\phi \mathcal{S} = - \mathcal{S} \phi$. Then, Lemma \ref{Hopf_thm_1} gives
		\begin{equation*}
			\phi \mathcal{S}^2 + c \phi = 0.
		\end{equation*}
		Therefore, for $\pi_\ast X \in U^\perp$ we have
		\begin{equation*}
			\phi (\mathcal{S}^2 + c I) \pi_\ast X = 0 = \phi(\lambda^2 + c)\pi_\ast X,
		\end{equation*}
		where $\lambda$ is eigenvalue of $\mathcal{S}$ corresponding to $\pi_\ast X$. Thus, we get
		\begin{equation*}
			\lambda^2 + c = 0.
		\end{equation*}
		For $c>0$, the preceding relation yields a contradiction, and hence $\sigma = 0$. When $c=0$, we obtain $\lambda=0$, in which case the shape operator $\mathcal{S}$ has precisely two eigenvalues, namely $\rho$ and $0$.
	\end{proof}
	
	\begin{remark}\label{remark1}
    Similarly to \cite{Ki_Suh}, we can prove that if $c < 0$, then $\mathcal{S} \phi + \phi \mathcal{S}$ cannot be identically zero. Hence, Lemma $\ref{thm_hopf_2}$ implies that $\rho$ is constant.
	\end{remark}
    
	Toward the nonexistence of umbilical hypermaps, we have the following theorem:
		\begin{theorem}\label{Nonexistence_umbilical_RHM}
		Let $\pi: (M, g_1) \to (N(4c), g_2, J)$ be a hypermap from a Riemannian manifold to a complex space form. If $c\neq 0$, then $\pi$ cannot be umbilical.
	\end{theorem}
	
	\begin{proof}
		Taking into account the almost contact structure $\{\phi, \eta\}$ and the Proposition \ref{umbilicity_cond_RM} in (\ref{Codazzi_RHM}), we obtain
		\begin{align*}
			&g_2( (\nabla_{\pi_\ast X}^N (g_1(Y, Z)H)), \xi) - g_2((\nabla_{\pi_\ast Y}^N (g_1(X, Z)H)), \xi) \\&= c \{ g_2(\phi \pi_\ast Y , \pi_\ast Z) \eta(\pi_\ast X)- g_2(\phi \pi_\ast X, \pi_\ast Z) \eta(\pi_\ast Y) + 2 g_2(\pi_\ast X, \phi \pi_\ast Y) \eta(\pi_\ast Z)\}.
		\end{align*}
		Since $\dim(\mathcal{R}^\perp) = 1$, the aforementioned equation gives
		\begin{align*}
			c \{ g_2(\phi \pi_\ast Y , \pi_\ast Z) \eta(\pi_\ast X)- g_2(\phi \pi_\ast X, \pi_\ast Z) \eta(\pi_\ast Y) + 2 g_2(\pi_\ast X, \phi \pi_\ast Y) \eta(\pi_\ast Z)\}=0.
		\end{align*}
		Now replacing $\pi_\ast Z$ with $\phi \pi_\ast X$, and then using (\ref{cc2}), we get
		\begin{align*}
			c \{ g_2(\pi_\ast X , \pi_\ast Y) \eta(\pi_\ast X) - g_2(\pi_\ast X , \pi_\ast X) \eta(\pi_\ast Y)\}=0.
		\end{align*}
		Taking $X \perp Y$, we get
		\begin{equation*}
			- c \| X \|^2 \eta(\pi_\ast Y)= 0,
		\end{equation*}
		which gives $c=0$. Thus, the required proof follows.
	\end{proof}

    Finally, we present the proof of Theorem $\ref{thm2_cpc_of_LR}$:
    
	\begin{proof}[Proof of Theorem $\ref{thm2_cpc_of_LR}$]
		Let $\mathcal{S}_\xi \phi + \phi \mathcal{S}_\xi = k \phi$ for some constant $k \neq 0$, then $U$ is an eigenvector of $\mathcal{S}_\xi$. Then (\ref{3.16}) yields,
		\begin{equation*}
			2 \phi \mathcal{S}_\xi^2 \pi_{\ast} X - 2k \phi \mathcal{S}_\xi \pi_{\ast}X + (\rho k + 2c) \phi \pi_{\ast}X= 0.
		\end{equation*} 
		On applying $\phi$, we get
		\begin{align*}
			2 \mathcal{S}_\xi^2 \pi_{\ast}X- 2 \mathcal{S}_\xi \pi_{\ast}X + (\rho k+2c) \pi_{\ast}X- 2(\rho^2 + c) \eta(\pi_{\ast}X) U + k \rho \eta(\pi_{\ast}X)U =0.
		\end{align*}
		Thus, the eigenvalues $\{\lambda_i\}_{i=1}^{2n-2}$ and $\rho$ of $\mathcal{S}_\xi$, satisfy the quadratic equation,
		\begin{equation*}
			2 t^2-2 kt +(\rho k+2c)=0.
		\end{equation*}
		Therefore, at most two eigenvalues $\lambda_i$ are distinct, and hence any $L_{\mathcal{R}}$ has at most three constant principal curvatures $\lambda, \mu$ (say) and $\rho$, as $\rho$, $k$, and $c$ are constant. If $\mathcal{S}_\xi \pi_{\ast} X = \lambda \pi_{\ast}X$, then
		\begin{equation*}
			\mathcal{S}_\xi \phi \pi_{\ast}X= (k-\lambda) \phi \pi_{\ast}X= \mu \phi \pi_{\ast}X.
		\end{equation*}
		Therefore, $\lambda$ and $\mu$ have equal multiplicities $n-1$. If $\lambda= \mu$, then
		\begin{equation*}
			\mathcal{S}_\xi \phi = \phi \mathcal{S}_\xi.
		\end{equation*}
		Hence, $2\mathcal{S}_\xi \phi= 2\phi\mathcal{S}_\xi= k \phi$. This implies that
		\begin{equation*}
			\mathcal{S}_\xi \pi_{\ast}X = \frac{1}{2}k \pi_{\ast}X + \frac{1}{2} (k-2 \rho) \eta (\pi_{\ast}X) U.
		\end{equation*}
		Hence, $\pi$ is $\eta$-umbilical hypermap.
	\end{proof}

    \section{Examples}
    In this section, we provide three examples of hypermaps to support the notions introduced in the current paper.

\begin{example}
Let $N =\mathbb{CP}^{2}(4c)$ be the complex projective space of real dimension $4$ with the complex structure $J$. By fixing a point $ p \in N$ and choosing $0 < t < \frac{\pi}{2\sqrt{c}}$, we have a real hypersurface $N_t$ in $N$, defined as
\[
N_{t} = \left\{ x\in N : r(x) = d(p,x)=t \right\},
\]
where $r(x)=d(p,x)$ denotes the geodesic distance from $p$ to $x$. In fact, $N_t$ is the geodesic sphere of radius $t$ centered at $p$, which is a Hopf and $\eta$-umbilical hypersurface in $N$.

We also define the set
\[
S^{3} = \left\{u\in T_{p} N : \|u\|=1\right\}.
\]

Equivalently,
\[
S^{3}
=
\left\{
(u_{1},u_{2},u_{3},u_{4})\in\mathbb{R}^{4}
:
u_{1}^{2}+u_{2}^{2}+u_{3}^{2}+u_{4}^{2}=1
\right\}.
\]

Then, one can see that the map $\pi: S^{3} \to N$ defined by 
\[
\pi(u) = \exp_{p} (t u)
\]
is a smooth map with $\pi (S^3) = N_t$. Furthermore, its derivative map $\pi_{*}: T_{u} S^{3} \to T_{\pi(u)} N$ satisfies $\pi_{*} (T_{u} S^{3}) = T_{\pi(u)} N_{t}$.
Thus, $\pi$ is a Hopf and an $\eta$-umbilical hypermap with 
$\mathcal{R}_x=T_x N_{t}$. 
One can also define a contact structure on $\mathcal{R}$ by taking the structure vector field $U$ as
\[U = -J (\nabla r).
\]
\end{example}

\begin{example}
    If we take $N = \mathbb{CH}^{2}(4c)$, a complex hyperbolic space, and choose $t>0$, then, similar to the previous example, the map $\pi: S^3 \to N$ also becomes a Hopf and an $\eta$-umbilical hypermap.
\end{example}






    \begin{example}
		Let $$\left( M =\{(x_1, x_2, x_3, x_4) \in \mathbb{R}^4 | -1 < x_1 < 1\}, ~g_1 = \sum \limits_{i=1}^{4} (dx_i)^2 + x_1 (dx_1) (dx_2)\right)$$ and $$\left(N =\{(y_1, y_2, y_3, y_4) \in \mathbb{R}^4 | -1 < y_1 < 1\}, ~g_2 = \sum \limits_{j=1}^{4} (dy_j)^2 + y_1 (dy_1) (dy_2)\right)$$ be two Riemannian manifolds. By \cite[{\rm ~Example} ~29]{Kiran_JMAA}, we have a Riemannian map $\pi : (M, g_1) \to (N, g_2)$ defined by 
            $$ \pi (x_1, x_2, x_3, x_4) = \left(x_1, 0, x_3, x_4\right)$$
            with 
            \begin{equation*}
				\mathcal{V} = \operatorname{span}\left\{e_2\right\}, \quad \mathcal{H} = \operatorname{span}\left\{e_1, e_3, e_4\right\},
			\end{equation*}
			\begin{equation*}
				\mathcal{R} = \operatorname{span}\left\{\pi_\ast (e_1) = e_1^*, \pi_\ast (e_3) = e_3^*, \pi_\ast (e_4) = e_4^*\right\}, \quad \mathcal{R}^\perp = \operatorname{span}\left\{e_2^*\right\},
			\end{equation*}
			where $\left\{e_i = \frac{\partial}{\partial x_i} \right\}_{i=1}^{4}$ and $\left\{e_j^* = \frac{\partial}{\partial y_j} \right\}_{j=1}^{4}$ are bases of $T_p M$ and $T_{\pi(p)} N$, respectively. Defining a complex structure $J$ on $N$ by $J(y_1, y_2, y_3, y_4) = (-y_3, -y_4, y_1, y_2)$, we have
		\begin{equation*}
			J e_2^* = -e_4^*.
		\end{equation*}
		Hence, $\pi$ is a hypermap with the structure vector field $U = e_4^*$. Further, we set $\xi = e_2^*$, $\eta( \pi_\ast X) = g_2(\pi_\ast X, e_4^*)$ and define $\phi: \mathcal{R} \to \mathcal{R}$ by $\phi \pi_\ast X = J \pi_\ast X - \eta(\pi_\ast X) e_2^*$. This yields
		\begin{equation*}
			\phi e_1^* = -e_3^*, \phi e_3^* = e_1^*, \phi e_4^* = 0, \eta(e_1^*) = 0, \eta(e_3^*) =0, \eta(e_4^*) = 1.
		\end{equation*}
		Thus, $\{\phi, \eta, U\}$ defines an almost contact structure.

        Furthermore, since $N$ is flat, one can easily verify that the map $\pi$ is a Hopf and an $\eta$-umbilical hypermap.
	\end{example}

    \section*{Acknowledgments}
	The authors gratefully acknowledge Dr. H. M. Shah for her participation in a few virtual discussions.
    
	
		
		
		
		

	\noindent K. Meena\\
	Department of Mathematics, Indian Institute of Technology Jodhpur, Rajasthan-342030, India.\\ E-mail: {kiranmeena@iitj.ac.in} 
    \\
    
    \noindent B. \c{S}ahin\\
	Department of Mathematics, Faculty of Science, Ege University, Izmir-35100, Turkiye.\\
	E-mail: {bayram.sahin@ege.edu.tr}
    \\

\addcontentsline{toc}{section}{References}
	\begin{thebibliography}{99}
		\bibitem{Akyol_2023} Akyol, M.A., G\"und\"uzalp, Y.: Pointwise slant Riemannian maps $(\mathcal{PSRM})$ to almost Hermitian manifolds. Mediterr. J. Math. {\bf 20} (2023), 1-15.
        
        \bibitem{Berndt} Berndt, J.: {Real hypersurfaces with constant principal curvatures in complex hyperbolic space}. J. Reine Angew. Math. \textbf{395} (1989), 132-141.

		

        \bibitem{Cecil_Ryan_Book} Cecil, T.E., Ryan, P.J.: Geometry of Hypersurfaces. Springer, New York, (2015).
		
		\bibitem{Chen_2001} Chen, B.Y., Maeda, S.: Hopf hypersurfaces with constant principal curvatures in complex projective or complex hyperbolic spaces. Tokyo J. Math. \textbf{24}(1) (2001), 133-152.
		
		\bibitem{Choe} Choe, Y.-W.: Characterizations of certain hypersurfaces of a complex space form. Nihonkai Math. J. \textbf{6} (1995), 97-114.


        \bibitem{Demir_2026} Demir, R., Kazan, S., Y\i ld\i r\i m, C.: On the structure of pointwise semi-slant Riemannian maps from Sasakian manifolds to Riemannian manifolds with vertical Reeb vector fields. Math. Nach. (2026). \url{https://doi.org/10.1002/mana.70280}
		
		\bibitem{Fischer_1992} Fischer, A.E.: {Riemannian maps between Riemannian manifolds}. Contemp. Math. \textbf{132} (1992), 331-366.
		
		\bibitem{Ki_Suh} Ki, U.-H., Suh, Y.J.: On real hypersurfaces of a complex space form. Math. J. Okayama Univ. \textbf{32} (1990), 207-221.	
		
		\bibitem{Kimura} Kimura, M.: Real hypersurfaces and complex submanifolds in complex projective space. Trans. Amer. Math. Soc. \textbf{296}(1) (1986), 137-149.
		
		\bibitem{Kimura_Maeda} Kimura, M., Maeda, S.: On real hypersurfaces of a complex projective space. Math. Z. \textbf{202}(3) (1989), 299-311.
		
		\bibitem{Kon} Kon, M.: {Pseudo-Einstein real hypersurfaces in complex space forms}. J. Differ. Geom. \textbf{14} (1979), 339-354.

        \bibitem{Lee_2021} Lee, C.W., Lee, J.W., \c{S}ahin, B., V\^{\i}lcu, G.-E.: Optimal inequalities for Riemannian maps and Riemannian submersions involving Casorati curvatures. Ann. Mat. Pura Appl. {\bf 200}(4) (2021), 1277-1295.
        
        \bibitem{Lee_2022} Lee, J.W., Lee, C.W., \c{S}ahin, B., V\^{\i}lcu, G.-E.: Chen-Ricci inequalities for Riemannian maps and their applications. Contemp. Math. \textbf{777} (2022), 137-152.

        \bibitem{MGSKRB_2026} Maini, R., Gupta, G., Kumar, R., Sachdeva, R., Rani, R., Bhatia, S.S.: Chen-Ricci and DDVV-type inequalities for anti-invariant Riemannian maps. Differential Geom. Appl. \textbf{104} (2026), 1-26.
        
        
        \bibitem{Meena_Sahin_Shah} Meena, K., \c{S}ahin, B., Shah, H.M.: Riemannian warped product maps. Result. Math. \textbf{79}(2) (2024), 1-30.
        
        \bibitem{Meena_2025}  Meena, K., \c{S}ahin, B., Shah, H.M.:  Chen's first inequality for Riemannian maps to complex space forms and $\delta $-invariants. Period. Math. Hungar. {\bf 91}(2) (2025), 588-603.

        \bibitem{Meena_Sahin_RHS} Meena, K., \c{S}ahin, B.: Riemannian hypersubmersions from complex space forms. Differential Geom. Appl. \textbf{104} (2026), 1-23.

        \bibitem{MMV} Meena, K., Meena, K.C., V\^ilcu, G.-E.: Geometry, classification, and rigidity of Clairaut Riemannian hypersubmersions. (Submitted).
		
		\bibitem{Montiel} Montiel, S.: Real hypersurfaces of a complex hyperbolic space. J. Math. Soc. Japan {\bf 37}(3) (1985), 515-535.
		
		\bibitem{Montiel_Romero} Montiel, S., Romero, A.: On some real hypersurfaces of a complex hyperbolic space. Geom. Dedicata \textbf{20}(2) (1986), 245-261.
		
		\bibitem{NR_1997} Niebergall, R., Ryan, P.J: {Real hypersurfaces in complex space forms}. Tight and Taut Submanifolds \textbf{32} (1997), 233-305.
		
		\bibitem{Nore_1987} Nore, T.: {Second fundamental form of a map}. Ann. Mat. Pura Appl. \textbf{146} (1986), 281-310.
		
		\bibitem{Okumura} Okumura, M.: On some real hypersurfaces of a complex projective space. Trans. Amer. Math. Soc. \textbf{212} (1975), 355-364.

        \bibitem{Penrose} Penrose, R.: Physical spacetime and non realizable CR-structure. Proc. Symp. Pure Math. \textbf{39} (1983), 401-422. 



        \bibitem{Polat_2025} Polat, G., Lee, J.W., \c{S}ahin, B.: Optimal inequalities involving Casorati curvatures along Riemannian maps and Riemannian submersions for Sasakian space forms. J. Geom. Phys. {\bf 210} (2025), 1-15.

        \bibitem{Polat_Meena} Polat, M., Meena, K.: Clairaut semi-invariant Riemannian maps to K\"ahler manifolds. Mediterr. J. Math. \textbf{21}(3) (2024), 1-19.

        \bibitem{Prasad_2025} Prasad, R., Kumar, S.: Point-wise semi-slant Riemannian maps from Riemannian manifolds into Kaehler manifolds. Filomat {\bf 39}(35) (2025), 12537-12554.
		
		\bibitem{Sahin_2025} \c{S}ahin, B.: Riemannian maps of CR-submanifolds of Kaehler manifolds. Math. Slovaca \textbf{75}(5) (2025), 1229-1240.
        
        \bibitem{Sahin_2010} \c{S}ahin, B.: {Invariant and anti-invariant Riemannian maps to K\"ahler manifolds}. Int. J. Geom. Methods Mod. Phys. \textbf{7}(3) (2010), 337-355.
		
		\bibitem{Sahin_2011} \c{S}ahin, B.: {Semi-invariant Riemannian maps to K\"ahler manifolds}. Int. J. Geom. Methods Mod. Phys. \textbf{8}(7) (2011), 1439-1454.
		
		
		\bibitem{Sahin_book} \c{S}ahin, B.: {Riemannian Submersions, Riemannian Maps in Hermitian Geometry, and Their Applications}. Elsevier, Academic Press, (2017).

        \bibitem{Kiran_JMAA} Singh, R., Meena, K., Meena, K.C.: {General Casorati inequalities and implications for Riemannian maps and Riemannian submersions}. J. Math. Anal. Appl. \textbf{558}(1) (2026), 1-31.

        \bibitem{SMM_GCFI} Singh, R., Meena, K., Meena, K.C.: {General Chen's first inequality and applications for Riemannian maps}. Indag. Math. (2026). \url{https://doi.org/10.1016/j.indag.2026.10.003}

        \bibitem{SMM_Filomat} Singh, R., Meena, K., Meena, K.C.: {General Chen-Ricci inequalities for Riemannian submersions and Riemannian maps}. Filomat (To appear).
		

		
		\bibitem{Takagi} Takagi, R.: On homogeneous real hypersurfaces in a complex projective space. Osaka J. Math. {\bf 10} (1973), 495-506.
		
		\bibitem{Tashiro_Tachibana} Tashiro, Y., Tachibana, S.-I.: On Fubinian and $C$-Fubinian manifolds. Kodai Math. Sem. Rep. \textbf{15}(3) (1963), 176-183.
		



        \bibitem{Yadav_Meena_2023} Yadav, A., Meena, K.: Riemannian maps whose base manifolds admit a Ricci soliton. Publ. Math. Debrecen \textbf{103}(1-2) (2021), 115-139.

        \bibitem{Jyoti_2024} Yadav, J., Shanker, G., Polat, M.: Clairaut slant Riemannian maps to K\"ahler manifolds. Int. J. Geom. Methods Mod. Phys. \textbf{21}(5) (2024), 1-19.
		
		\bibitem{Yano} Yano, K., Kon, M.: Structure on Manifolds. World Scientific, Singapore, (1984).

	\end{thebibliography}
\end{document}